\documentclass[letterpaper, 10 pt, conference]{ieeeconf}   

\IEEEoverridecommandlockouts                              
\usepackage{graphicx}
\usepackage{subfigure}
\usepackage{amssymb}
\usepackage{widetext}
\usepackage{multicol}
\usepackage{amsmath} 
\usepackage{lipsum}  
\usepackage{cite}

\usepackage{color}

\title{\LARGE \bf
Identifiability of dynamical networks: which nodes need be measured? }

\author{Alexandre S.  Bazanella$^{1}$, Michel Gevers$^{2}$,  Julien M. Hendrickx$^{2}$ and Adriane Parraga$^3$
\thanks{*This work is supported
 by the Program Science Without Borders, CNPq - Conselho Nacional de Desenvolvimento Cient\'{\i}fico e Tecnol\'{o}gico, Brazil,   by the Belgian Programme on
Interuniversity Attraction Poles, initiated by the Belgian Federal
Science Policy Office and by Wallonie Bruxelles International.}
\thanks{$^{1}$A. S. Bazanella is with the Department of Automation and Energy, 
        Universidade Federal do Rio Grande do Sul, Brazil
        {\tt\small bazanella@ufrgs.br }}%
\thanks{$^{2}$M. Gevers is with ICTEAM, Univsersit\'e catholique de Louvain, B1348 Louvain la Neuve, Belgium        {\tt\small Michel.Gevers@uclouvain.be}}%
\thanks{$^{2}$J. Hendrickx is with ICTEAM, Univsersit\'e catholique de Louvain, B1348 Louvain la Neuve, Belgium        {\tt\small Julien.Hendrickx@uclouvain.be}}%
\thanks{$^{3}$A. Parraga  is with the Computer Engineering Department, Universidade Estadual do Rio Grande do Sul, Brazil
        {\tt\small adriane-parraga@uergs.edu.br }}%
}

\newcommand{\bydef}{\stackrel{\Delta}{=}}

\newcommand{\beq}{\begin{equation}}
\newcommand{\eeq}{\end{equation}}
\newcommand{\beqa}{\begin{eqnarray}}
\newcommand{\eeqa}{\end{eqnarray}}
\newcommand{\beqan}{\begin{eqnarray*}}
\newcommand{\eeqan}{\end{eqnarray*}}
\newcommand{\bef}{\begin{figure}}
\newcommand{\enf}{\end{figure}}

\newcommand{\bi}{\begin{itemize}}
\newcommand{\ei}{\end{itemize}}
\newcommand{\bc}{\begin{center}}
\newcommand{\ec}{\end{center}}
\newcommand{\ba}{\begin{array}}
\newcommand{\ea}{\end{array}}
\newcommand{\be}{\begin{equation}}
\newcommand{\ee}{\end{equation}}
\newcommand{\beno}{\begin{equation*}}
\newcommand{\eeno}{\end{equation*}}
\newcommand{\beqna}{\begin{eqnarray}}
\newcommand{\eeqna}{\end{eqnarray}}
\newcommand{\bd}{\begin{displaymath}}
\newcommand{\ed}{\end{displaymath}}
\newcommand{\beqnd}{\begin{eqnarray*}}
\newcommand{\eeqnd}{\end{eqnarray*}}

\newcommand{\ta}{\theta}

\renewcommand{\ni}{\noindent}

\newcommand{\cqfd}{\hfill \rule{1.8mm}{1.8mm}\medbreak\indent}

\newtheorem{theorem}{\bf Theorem}[section]
\newtheorem{lemma}{\bf Lemma}[section]

\newtheorem{corollary}{\bf Corollary}[section]

\usepackage{color}

\definecolor{red}{rgb}{1,0,0}

\definecolor{blu}{rgb}{0,0,1}

\definecolor{gre}{rgb}{0,0.7,0.3}

\begin{document}

\maketitle

\begin{abstract}
Much recent research has dealt with the identifiability of a dynamical network in which the node signals are connected by causal linear time-invariant transfer functions and are possibly excited by known external excitation signals and/or unknown noise signals. So far all  results on the identifiability of the whole network have assumed that all node signals are measured. Under this assumption, it has been shown that such networks are  identifiable only if some prior knowledge is available about the structure of the network, in particular the structure of the excitation. In this paper we present the first results for the situation where not all node signals are measurable, under the assumptions that the topology of the network is known, that each node is excited by a known signal and that the nodes are noise-free. Using graph theoretical properties, we show that the transfer functions that can be identified depend essentially on the topology of the paths linking the corresponding vertices to the measured nodes. An important  outcome of our research is that, under those assumptions, a network can often be identified using only a small subset of node measurements. 
\end{abstract}

\section{INTRODUCTION}\label{intro}
This paper examines the identifiability of dynamical networks in which the node signals  are connected by causal linear time-invariant transfer functions and are possibly excited by known external excitation signals and/or unknown noise signals. Such networks can be looked upon as directed graphs in which the edges between the nodes (or vertices) are composed of scalar causal linear transfer functions, and in which known external excitation signals enter into the nodes.

The identification of  networks of linear time-invariant dynamical systems based on the measurement of all its node signals and of all known external excitation signals acting on the nodes has been the subject of much recent attention \cite{Goncalves&Warnick:08,Materassi&Innocenti:10,Dankers&Vandenhof&Heuberger&Bombois:12,Chiuso&Pillonetto:12, Weerts&Dankers&Vandenhof:15,Hayden&Chang&Goncalves&Tomlin:16,Gevers&Bazanella&Parraga:17}. It has been shown in \cite{Goncalves&Warnick:08,Gevers&Bazanella&Parraga:17} that, generically, such networks cannot be identified from the node signals and the known external excitation signals, and that identifiability can only be obtained provided prior knowledge is available about the structure of the network. In practice, it is often the case that the  excitation structure is known, i.e. one  often knows at which nodes  external excitation signals or unknown noise signals are applied. A number of sufficient conditions for the identifiability of the whole network have been derived under prior assumptions on the structure of the network, involving either its external excitation structure,  or possibly also its internal structure \cite{Goncalves&Warnick:08,Weerts&Dankers&Vandenhof:15,Hayden&Chang&Goncalves&Tomlin:16,Gevers&Bazanella&Parraga:17}.

In all the results accumulated so far on the identifiability of networks of dynamical systems, it is assumed that all node signals can be measured. In this paper we examine the situation where not all node signals are measured, but where the topology of the network is known; this means that the user knows a priori which nodes are connected by nonzero transfer functions. We also assume that  each node is excited by a known external excitation signal, but that the node signal is noise-free.
In this context, a number of questions can be raised, such as
\begin{enumerate}
\item Can one identify the whole network with a restricted number of node measurements?
\item If so,  are there a minimal  number of nodes that need to be measured?
\item Are some node measurements useless, in  that they deliver no information about any transfer function?
\item Are some nodes indispensable, in  that it is impossible to identify the network without measuring them?
\item If one wants to identify a specific transfer function, can the topology tell us which node or nodes need to be measured?
\item Which transfer functions can be identified from the measure of a specific subset of nodes?
\end{enumerate}

To answer these questions we shall heavily rely on properties from graph theory, using the directed graph corresponding to our network as our major tool. 

To the best of our knowledge, the only other contributions that consider identification in networks with unobservable nodes are \cite{Materassi&Salapaka:15,Dankers&Vandenhof&Bombois&Heuberger:16,Linder&Enqvist:17}.
However, the problem treated in these papers consists of 
the identification of a subset of the networkÕs transfer functions - typically a single one - and hence
is different from the one presented here. 
Sparse measurements have also been considered in a different context in \cite{Mauroy&Hendrickx:16};
the goal there was to recover the network structure under the assumptions that the local dynamics are known,
as opposed to re-identifying the dynamics and/or the whole network structure.

After stating the problem in precise terms in Section~\ref{probstatement}, we shall first, in Section~\ref{motivation}, motivate the reason for addressing this problem by exhibiting some simple 3-node networks to show that the nodes that need to be measured to identify all transfer functions depend completely on the topology of the network and that, in some cases, a unique measurement suffices to identify the whole network for a 3-node network. This already yields a positive answer to question 1 above. Our brief analysis of 3-node networks will then lead us, in Section~\ref{basicresults} to formulate a number of basic results pertaining to questions 2, 3 and 4 above.

In Section~\ref{pathbased} we shall present a number of identifiability results that are based on properties of the paths going from a specific node within the network to a specific set of measured nodes. This will lead to a necessary and sufficient result on identifiability of all transfer functions, or of a specific set of transfer functions, that are based on the properties of the path from each node to a measured node. This section will also give a partial answer to question 5 above.

In Section~\ref{measurementbased} we shall address question 6 above. Instead of looking at a specific node within the network and examining its paths to a measured node or a set of measured nodes, as was done in Section~\ref{pathbased}, we consider the converse approach. We look at a specific measured node and we present a  graph-theoretical result that lets us decide which transfer functions can be identified from that specific node.

In Section~\ref{conclusion} we will conclude and describe  some challenging  open problems that remain to be solved.

\section{Statement of the problem}\label{probstatement}
The problem studied in this paper is part of the recent research on the question of identifiability of networks of dynamical systems. We first present the network structure and explain the network identifiability problem as it has so far been posed, i.e.  with all nodes measured. We then focus on the new network identifiability problem for the case when not all nodes are measured under the simplifying assumption that there is a known external excitation signal that acts on each node.

We adopt a simplified noise-free version of the standard network structure of  \cite{Weerts&Dankers&Vandenhof:15,Gevers&Bazanella&Parraga:17} for networks whose edges are linear causal scalar rational transfer functions. 
Thus, we consider that the network is made up of $L$ nodes, with  node signals  denoted $\{w_1(t), \ldots, w_L(t)\}$, and that  these node signals are related to each other and to  external excitation signals $r_j, j=1,\dots,L$ by the following network equations, which we call the {\bf network model} and in which the matrix $G^0$ will be called the {\bf  network matrix}:
\beqnd 
\hspace{-2mm}\left[ \begin{array}{c}
w_1\\ w_2\\ \vdots \\ w_L \end{array} \right] \!=\! \left[\! \begin{array}{cccc}0 & \!\! G_{12} & \!\! \ldots &\!\!  G_{1L} \\G_{21} & \!\! 0  \!\!  & \ddots & \!\! G_{2L} \\Ê\vdots &\!\! \ddots &\!\!  \ddots &\!\!  \vdots\\
G_{L1} & \!\! G_{L2} & \!\! \ldots & \!\! 0\end{array}\! \right] \!\! \left[\!  \begin{array}{c}
w_1\\ w_2\\ \vdots \\ w_L \end{array} \! \right] \!+\! K^0(q) \left[ \begin{array}{c}
r_1\\ r_2\\ \vdots \\ r_L \end{array} \right] 
\eeqnd
Or, equivalently
\be
w(t) = G^0(q) w(t) + K^0(q) r(t)  \label{netmodel} 
\ee
with the following  properties.
\begin{itemize}
\item $G_{ij}$ are proper but not necessarily strictly proper transfer functions. Some of them may be zero, indicating that there is no direct link from $w_j$ to $w_i$.
\item there is a delay in every loop going from one $w_j$ to itself.
\item the network is well-posed so that $(I-G^0)^{-1}$ is proper and stable.
\item $r_i$ are  known external excitation signals that are available to the user in order to produce informative experiments for the identification of  the $G_{ij}$. They are asumed to be mutually uncorrelated. The $L \times L$ transfer function matrix $K^0(q)$ reflects how the external excitation signals affect the node signals. 
\item $q^{-1}$ is the delay operator.
\item the topology of the network is  known, i.e. one knows which of the $G_{ij}$ are zero.
\end{itemize}

The network model (\ref{netmodel}) can be rewritten in a more traditional form as follows:
 \be 
w(t) = T^0(q) r(t)  \label{iomodel2}
\ee
where 
\be \label{TNdef}
T^0(q) \bydef (I - G^0(q))^{-1}K^0(q).
\ee
The description (\ref{iomodel2}) will be called the {\bf input-output (I/O) description} of the network. 
A corresponding parametrized version $M_{io} = T(q,\ta)$ will be called the {\it input-output (I/O) model}. 

The question that has attracted a lot of attention in the last few years is that of the identifiability of the network and of the informativity of the external excitation data $r(t)$. The question of network identifiability can  be briefly stated as follows. Assuming that the network is driven by sufficiently informative excitation signals $r(t)$, what are the conditions (in the form of required prior knowledge) on the network matrices $G^0(q), K^0(q)$ such that they can be uniquely identified from measurements of the node signals $w(t)$ and the known external signals $r(t)$? Assuming that the network structure is identifiable, the question of informativity is then: what excitation is required from the external signals $r(t)$ for a parametric and identifiable network model structure $[G(q,\ta), K(q,\ta)]$ to converge to the true  $[G^0(q), K^0(q)]$.

It is well known from the theory of identification of multi-input multi-output (MIMO) linear time-invariant (LTI) systems that on the basis of measurements of the signals $w(t)$ and $r(t)$ one can uniquely identify the matrix $T^0(q)$  of the input-output model (\ref{iomodel2}) if the chosen model structure $M_{io} = T(q,\ta)$ is identifiable and such that $T^0(q)= T(q,\ta_0)$ for some $\ta_0$  (this is the identifiability question), and if  the signals $r(t)$ are sufficiently rich for the chosen parametrizations (this is the informativity question).
Observe that the identification of (\ref{iomodel2}) is an open loop identification problem for a LTI MIMO system. 

The question of {\bf network identifiability} then relates to the mapping from $T^0(q)$ to $[G^0(q), K^0(q)]$, namely {\it under what conditions (in the form of prior knowledge on the network matrices 
$G^0(q), K^0(q)$) can one uniquely recover $[G^0(q), K^0(q)]$ from $T^0(q)$?}

So far, all results on this problem have been obtained under the assumptions that all nodes $w_i(t)$ are measured without error. A range of sufficient conditions have been obtained for this situation: see \cite{Goncalves&Warnick:08, Weerts&Dankers&Vandenhof:15,Gevers&Bazanella&Parraga:17}. A representative example of such sufficient conditions states that the network (\ref{netmodel}) is identifiable if  $K^0(q)$  is diagonal and of full rank.

In this paper we address the question of network identifiability for the case where not all nodes are measured, but where the topology of the network is known. We consider that, together with the network (\ref{netmodel}), there is a measurement equation
\be \label{measures}
y(t) = C w(t)
\ee
where $C$ is a $pÊ\times L$ matrix that reflects the selection of measured nodes. Thus, each row of $C$ contains one element $1$ and $L-1$ elements $0$. We shall denote by $\cal C$ the corresponding subset of nodes selected by $C$.  Given that the problem turns out to be difficult, we start in this paper with the simple situation where $K^0(q) = I_L$,  and the excitation signals  $r_i(t)$ are assumed to be independent of one another and sufficiently rich of any desired degree. Under those assumptions, the network simplifies to 
\beqna 
w(t) & = & G^0(q) w(t) + r(t) \label{netmodel3} \\
y(t) &=& Cw(t) = CT^0(q) r(t) \label{netmodel2}
\eeqna
where
\be
T^0(q) = (I - G^0(q))^{-1}. \label{Tdef}
\ee
We thus summarize the assumptions that will be made throughout this paper.\\
\ni {\bf Standing assumptions.}
\begin{itemize}
\item The networks  and the measurements are described by (\ref{netmodel3})-(\ref{netmodel2}), where $G^0(q)$ has the properties described above.
\item The  topology is known, i.e. one knows which of the $G_{ij}(q)$ are zero.
\item A sufficiently rich external excitation  signal $r_i(t)$ is applied to each node, and these signals are mutually uncorrelated.
\end{itemize}

Clearly, the matrix $CT^0(q)$ can be uniquely identified from $\{y(t), r(t)\}$ data. The network identifiability problems addressed in this paper are then as follows: {\it Which transfer functions $G^0_{ij}(q)$ can we uniquely identify from $CT^0(q)$, and under what conditions can we identify the whole network matrix $G^0(q)$ from $CT^0(q)$?}  

The identification of the $G_{ij}(q)$ rests on the following relationship
\be \label{TzeroT1}
C T^0(q) = CT(q) = C(I-G(q))^{-1}
\ee
or, equivalently,
\be \label{TzeroT2}
C T^0(q) (I-G(q)) = C
\ee
For a given $C$, $C T^0(q)$ is assumed known since it can be perfectly identified from $\{y,r\}$ data. We then  solve (\ref{TzeroT2})   for the unknown $G_{ij}$ and check which of these can be identified. We note from (\ref{TzeroT2}) that each node measure contributes $L$ scalar equations in the scalar unknowns $G_{ij}$. However, as we shall see below, depending on the network topology, some of these equations may be trivial equations (such as 1=1, or 0=0) and therefore do not contain any information about the $G_{ij}$.

\section{Motivating examples}\label{motivation}
In order to motivate the reader, we now analyze a few 3-node networks and show that the nodes that allow identification of the whole network depend entirely on the topology of the network, and that the entire network can often be identified from the measurements of a small subset of nodes. 

\begin{figure}[h]
 \centering
 \begin{tabular}{cc}
 \begin{minipage}{.4\textwidth}
\includegraphics[width=0.4\linewidth]{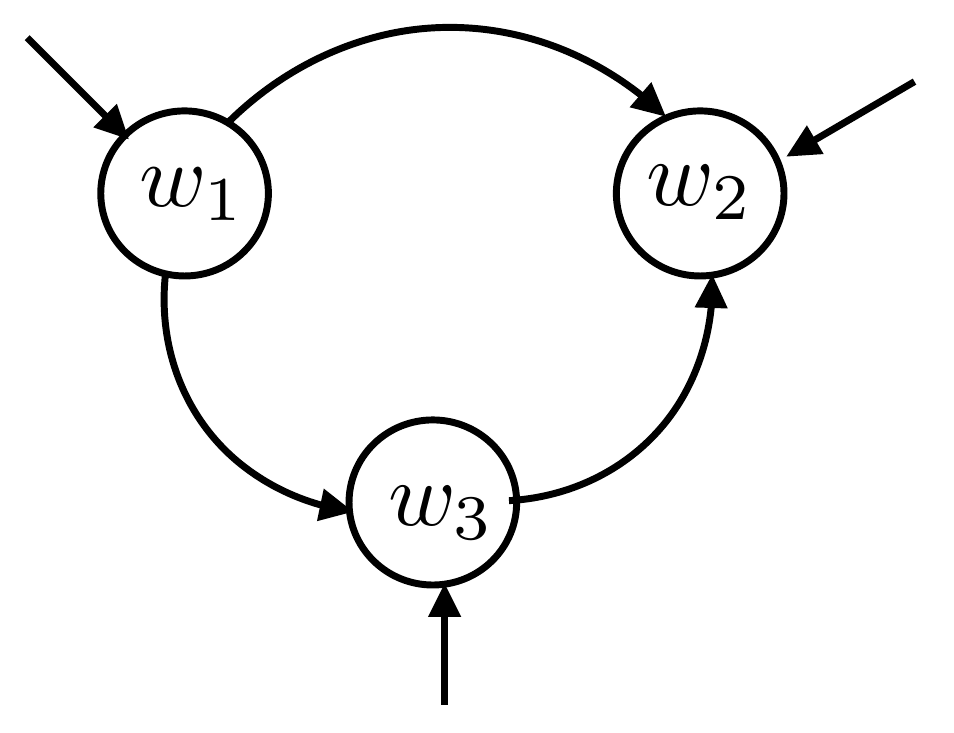} 
    \end{minipage}
  $\hspace{-40mm}  G(q) =\left[\begin{array}{ccc}0 & 0 & 0 \\ G_{21} & 0 & G_{23}\\ G_{31} & 0 & 0 \end{array}\right], $ \\     
$T(q) =  \left[\begin{array}{ccc}1 &0 & 0 \\ G_{21} +  G_{23}G_{31}& 1 & G_{23}\\ G_{31}&0 &1 \end{array}\right] $
\end{tabular}
\caption{Example of network with three transfer functions where two nodes (2 and 3) need to be measured.}\label{fig:network1}
\end{figure}

Consider first a network with 3 unknown transfer functions represented in Figure \ref{fig:network1} and its corresponding true $G(q)$ and true $T(q)$. Calculations based on  (\ref{TzeroT2}) show that identification of all 3 transfer functions requires the measurement of nodes 2 AND 3, and that measuring node 1 yields no information.

\begin{figure}[h]
 \centering
 \begin{tabular}{cc}
   \begin{minipage}{.4\textwidth}
\includegraphics[width=0.4\linewidth]{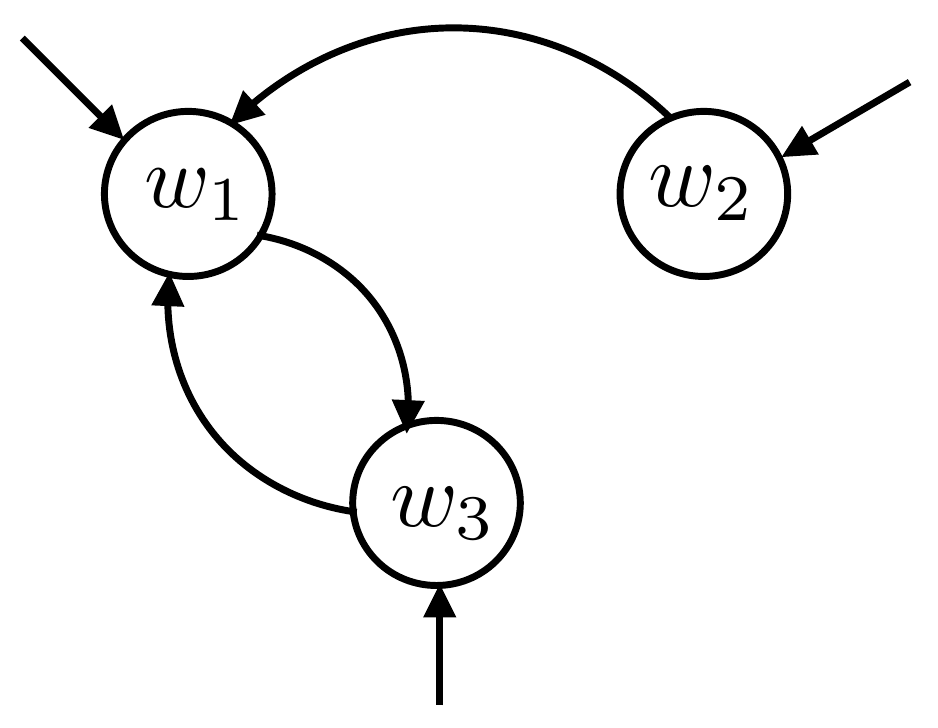} 
    \end{minipage}
 $\hspace{-40mm} G(q) =\left[\begin{array}{ccc}0 & G_{12} & G_{13} \\ 0 & 0 & 0 \\ G_{31} & 0 & 0 \end{array}\right], $ \\
 $T(q) =  \frac{1}{\Delta} \left[\begin{array}{ccc}1 &G_{12}  & G_{13} \\ 0 & 1- G_{13}G_{31} &0 \\ G_{31}  & G_{31}G_{12} & 1 \end{array}\right] $
\end{tabular}
\caption{Example of network with three transfer functions where measuring one node (1 or 3) is sufficient. We use $\Delta \bydef det(I-G)= 1 - G_{13}G_{31}$}\label{fig:network2}
\end{figure}

By contrast, the identification of the 3 unknown transfer functions in the network represented in Figure \ref{fig:network2}  is possible by measuring just one node: node 1 OR node 3.

\begin{figure}[h]
 \centering
 \begin{tabular}{cc}
 \begin{minipage}{.4\textwidth}
\includegraphics[width=0.4\linewidth]{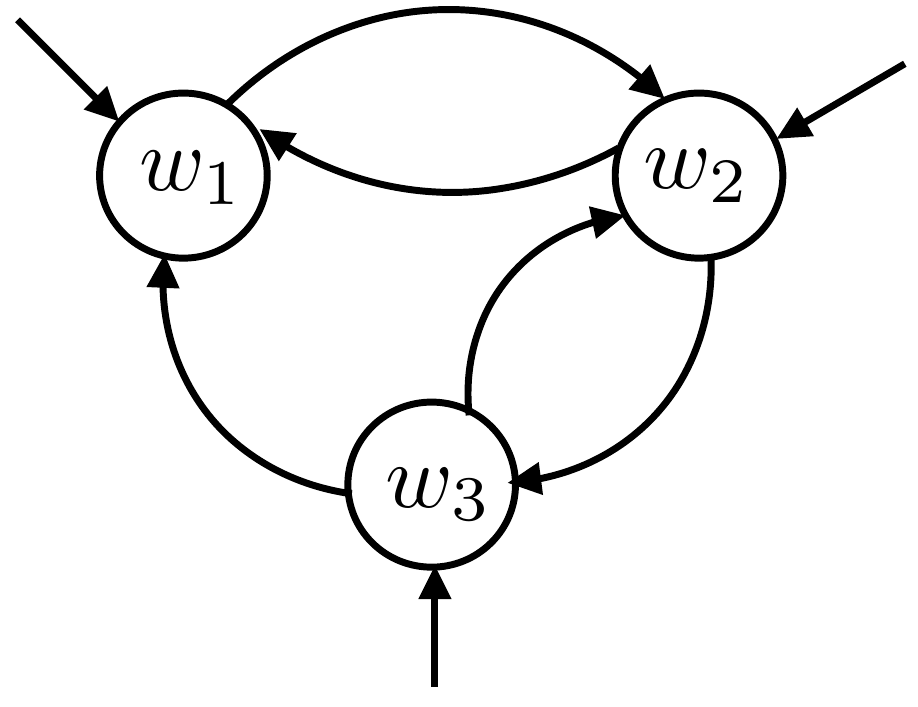} 
    \end{minipage}
 $\hspace{-40mm} G(q) =\left[\begin{array}{ccc}0 & G_{12} & G_{13} \\ G_{21} & 0 & G_{23} \\ 0 & G_{32} & 0 \end{array}\right], $ \\     
\hspace{-3mm} $T(q) \!\!=\!\!  \frac{1}{\Delta}\! \left[\begin{array}{ccc}1-G_{23}G_{32}  & G_{12}+G_{13}G_{32}    & G_{13}+ G_{12}G_{23} \\ G_{21} & 1 & G_{23} + G_{21}G_{13}  \\ G_{32}G_{21}  & G_{32}  & 1-G_{12}G_{21}  \end{array}\!\! \right] $
\end{tabular}
\caption{Example of network with five transfer functions where measuring two nodes (1 and 2 or 1 and 3) is sufficient. We use $\Delta \bydef  1-G_{12}G_{21} -  G_{23}G_{32}- G_{13}G_{21}G_{32}  $.}\label{fig:network3}
\end{figure}

Finally,  in the network of Figure \ref{fig:network3}, all  5 transfer functions can be identified by measuring  just two nodes:  either nodes 1 AND 2, OR nodes 1 AND 3.  

The examples show  that the number of measurements that are necessary to identify the network depends not only on the number of unknown
transfer functions to be determined (the number of nonzero $G_{ij}$) but also on the topology of the network; it is also clear that not all measurements are equal.

\section{Basic results}\label{basicresults}
Inspired by our analysis of 3-node networks, we now establish a number of basic results regarding the identifability of general $L$-node networks from a reduced set of node measurements. 

 We first introduce some  notations and definitions from graph theory (see e.g. \cite{Diestel:00}). \\
{\bf Notations and definitions:} 
\begin{itemize}
\item $L$ = number of nodes;  
\item $p$ = number of measured nodes; 
\item $f$ = number of sources, i.e. number of nodes with only outgoing edges;  
\item $s$ = number of sinks, i.e. number of nodes with only incoming edges;  
\item $n$ = number of unknown transfer functions; 
\item $C$ = the $p \times L$ matrix that reflects the selection of nodes via $y(t) = C w(t)$: thus  each row of $C$ contains one element $1$ and $L-1$ elements $0$; 
\item $\cal C$ = the subset of nodes selected by $C$;
\item $N_i^+$ = set of out-neighbors  of node $i$, i.e. the set of nodes $j$ for which $G_{ji}\neq 0$;
\item  $d_i^+$ = number of outgoing edges of node $i$, i.e. number of nodes in $N_i^+$;
\item  a {\it walk}Ê~denotes a series of adjacent directed edges (including trivial walks consisting of one node with no edge);
\item a {\it cycle} is a walk whose terminal node coincides with the initial node;\footnote{Note that a cycle is typically called {\it loop} in control theory, but we shall from here on use cycle, which is standard in graph theory.} 
\item  a {\it path} is a 
walk that never passes twice through the same node, i.e. a walk without loops. 
\end{itemize}

We can now establish the  following basic results. 
\begin{theorem}\label{basic1}
1) If $w_i$ is a source, then $G_{ij} =0 \;\forall j$, $T_{ii}=1$ and $T_{ij}=0 \;\forall j\neq i$.  It then follows from (\ref{TzeroT2})  that the corresponding equation yields no information. Thus, measurement of sources is useless, but each source requires excitation.\\
2) If $w_i$ is a sink, then $G_{ji} =0 \;\forall j$, $T_{ii}=1$ and $T_{ji}=0 \;\forall j\neq i$. Identifiability of the network requires that all sinks be measured, but exciting a sink is useless.\\
\ni{Proof:} see Appendix.
\end{theorem}

We now make some observations concerning the number of useful equations that result from (\ref{TzeroT2})  for the computation of the $G_{ij}$. Each measured node contributes $L$ equations, but some of these may not yield any information.

First we note that $L-1 \leq n \leq L(L-1)$. The number of equations is $p\times L$; so it is obvious
that we need $p \geq \frac{n}{L}$. It now follows from (\ref{TzeroT2}) and Theorem~\ref{basic1} that each  sink causes the appearance of one trivial equation $1=1$ in the sink's measurement, and also 
of one trivial equation $0=0$ at every other measurement. Hence the number of trivial equations 
caused by each sink equals $p$, and thus the total number of trivial equations
due to the existence of sinks is $sp$. Therefore  the number of useful equations is
at most $n_e = pL -  ps = p(L-s)$. We then have the following result.

\begin{theorem}\label{basic2}
Identifiability of the whole network requires measurement of all sinks plus at least $m$ more nodes such that 
\be \label{minnumber}
m+s \geq \frac{n}{L-s}
\ee
\ni{\bf Proof:} Given that the number of useful equations resulting from $p$ measurements is $p(L-s)$, identifiability of a  network with $n$ unknowns and $s$ sinks requires that $p(L-s)\geq n$, where $p=m+s$. This implies (\ref{minnumber}).
\end{theorem}
The next theorem yields a simple result for a  class of acyclic networks that contains the trees as important particular case (These are networks that cannot be separated in disjoint parts, and do not contain any cycle even if one were to change the direction of some edges).
\begin{theorem}\label{treetheo}
Consider a  network that contains no cycle, and where no node is connected to any other one by more than one path. Then measuring all sinks is necessary and sufficient for the identifiability of all unknown transfer functions.
\end{theorem}
The next result covers the case of loop structures.
\begin{theorem}\label{theoloop}
Let the nodes $w_i$, $i \in {\cal I}$  form one cycle and assume that 
 no other cycle in the graph contains any of these nodes. Then measuring any one of these nodes is sufficient to identify all transfer functions in the cycle. \\
 \ni{Proof:} see Appendix.
 \end{theorem}

\section{Path-based results}\label{pathbased}
In this section we present a series of identifiability results that are based on the structure of the paths from a given node to a set of measured nodes. We first reformulate the identifiability problem. Recall that $G$ is identifiable if and only if (\ref{TzeroT2}) implies $G=G^0$ for any $G$ \emph{consistent with the graph, i.e. with the topology}. Define $\Delta \bydef G-G^0$, which is consistent with the graph if and only if $G$ is. Substituting $G=G^0 + \Delta$ in (\ref{TzeroT2}) shows that $G$ is identifiable if and only if
\be \label{eq:reformulation_Delta}
C T^0(q) \Delta(q) =0 \Rightarrow \Delta(q) = 0
\ee
for any $\Delta(q)$ consistent with the graph. 

Condition \eqref{eq:reformulation_Delta} can be viewed as a set of $L$ conditions, one for each column, which must all be true. The condition for the $i$-th column is: 
\begin{equation} \label{eq:condition_columnn}
C(I-G^0)^{-1}\Delta_{:i} = 0 \Rightarrow \Delta_{:i} = 0,
\end{equation}
for $\Delta_{:i}$ consistent with the graph i.e. $\Delta_{ki} = 0$ if there is no edge $(i,k)\in G$.  Condition \eqref{eq:condition_columnn} can be rewritten as
\begin{equation}\label{eq:condition_columnn_detail}
\sum_{k\in N^+_i} T^0_{jk}\Delta_{ki} = 0, \forall j\in {\cal C}, \Rightarrow \Delta_{ki} = 0, \forall k\in N^+_i
\end{equation}
Note that one can prove that $T^0_{jk}$ is nonzero only if there is a path from $k$ to $j$.
Note also that if condition \eqref{eq:condition_columnn_detail} is satisfied for some $i$, then we know that $G_{:i} = G^0_{:i}$ for any $G$ consistent with the graph and satisfying $C(I-G)^{-1} = C(I-G^0)^{-1}$, i.e. we can identify all transfer functions corresponding to the edges leaving node $i$.

We now establish a number of identifiability results based on the paths from a given node $i$ to a set of measured nodes.

\begin{theorem}\label{prop:suf_condition}
For a node $i$, if each of  its $d^+_i$ out-neighbors has a directed path to a different measured node $j\in {\cal C}$, and if these directed paths are all vertex-disjoint, including their measured end-points,   then condition \eqref{eq:condition_columnn_detail} is generically satisfied (i.e. it is satisfied for almost all choices of transfer functions matrices $G^0$ consistent with the graph), and the transfer functions corresponding to the edges leaving $i$ can all be identified. \\
\ni{Proof:} see Appendix.
\end{theorem}

We stress that the sufficient condition in Theorem \ref{prop:suf_condition} does not require all paths from the $d^+_i$ out-neighbors of $i$ to measured nodes to be disjoint, but only the existence of a set of mutually disjoint paths. In other words, there may very well exist many other paths than those used in the condition, and there is no requirement on those, nor on their intersections with those used in the condition. 

The next result shows sufficient conditions under which a subset of the out-going edges can be recovered even if the others cannot. In particular, it yields a solution for the identifiability of a single embedded transfer function.
\begin{theorem}\label{prop:suf_condition_partial}
Consider a node $i$, and let $N^*_i\subseteq N^+_i$ be a subset of its out-neighbors. Suppose in addition that

(i) There exists $d^*_i = |N_i^*|$ vertex-disjoint directed paths joining nodes of $N^*_i$ to $d^*_i$ measured nodes, and let ${\cal C^*}$ be the set of these measured nodes.  

(ii) There is no path from any node of $N^+_i\setminus N^*_i$ to any node of ${\cal C^*}$.

\ni Then the transfer functions corresponding to edges from $i$ to nodes in $N^*$ can be identified.\\
\ni{Proof:} see Appendix.
\end{theorem} 

We now present a necessary condition. 
\begin{theorem}\label{prop:nec_condition}
Consider a node $i$ and suppose there exists a set of nodes $B$  with cardinality $b<d^+_i$ such that any path from an out-neighbor $i\in N^+_i$  of $i$ to a measured node $j\in {\cal C}$ includes a node in $B$.\footnote{ $B$ stands for bottleneck; note that it may contain nodes of $N^+_i$.}
Then condition \eqref{eq:condition_columnn_detail} is generically not satisfied, and it is therefore not possible to identify the transfer functions of all edges leaving $i$. \\
\ni{Proof:} see Appendix.
\end{theorem}

The absence of such a set $B$ is thus a necessary condition for identifiability of all outgoing edges of node $i$.

The following Lemma implies that that there is no gap between the sufficient condition of Theorem~\ref{prop:suf_condition} and the necessary condition of Theorem~\ref{prop:nec_condition}.

\begin{lemma}\label{lem:duality}
Let $S$ and $P$ be sets of nodes in a directed graph, possibly with a non-empty intersection. For any $d$, exactly one of the two following conditions holds:\\
(i) There exist $d$ vertex-disjoint paths, each joining a node of $S$ to a node of $P$.\\
(ii) There is a set $B$ of $b< d$ nodes such that every path from $S$ to $P$ contains a node of $B$.
\end{lemma}
\ni {\bf Proof:}  This follows from a variation of the min-cut max-flow duality \cite{Papadimitriou&Steiglitz:98}. \cqfd

\begin{figure}[h]
 \centering
\includegraphics[width=0.7\linewidth]{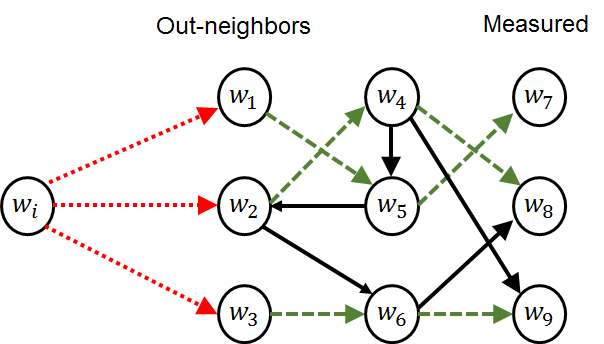}     
\caption{Example illustrating Theorem~\ref{thm:nec_suf_condition}: 3 vertex-disjoint dashed green paths to the 3 measured nodes; the 3 dotted red edges are identifiable.}
\label{fig:nscillustration}
\end{figure}

Our main path-based result provides necessary and sufficient conditions for the identifiability of all outgoing edges of any node $i$ in a network.
\begin{theorem}\label{thm:nec_suf_condition}
For every node $i$, one can generically uniquely identify all transfer functions on the edges leaving $i$ if and only if there exist vertex-disjoint directed paths leaving the out-neighbors of $i$ and arriving at measured nodes. \footnote{The vertex-disjoint condition applies also for the departure and arrival nodes.}
In particular, (\ref{TzeroT2}) yields $G(q)=G^0(q)$  if and only if the condition above is satisfied for every $i$.
\end{theorem}
\ni {\bf Proof:} This follows from Lemma \ref{lem:duality} and Theorems  \ref{prop:suf_condition} and \ref{prop:nec_condition}. \cqfd
The theorem is illustrated by the  example in Figure~\ref{fig:nscillustration}. Remember that known external signals $r_i$ are applied to each node, which we have not added on the figure for visibility reasons. Node $i$  has three outgoing nodes, each of which has a vertex-disjoint directed path to the measured nodes 7, 8 and 9, namely the paths $(1,5,7), (2,4,8)$ and $(3,6,9)$; they are represented by dashed green arrows. As a result, the dotted red transfer functions $G_{1i}, G_{2i}$ and $G_{3i}$ can all be identified from these three measured  nodes.

An immediate consequence is the following necessary condition for identifiability of the whole network.
\begin{corollary}\label{coro1}
The number of measured nodes must be larger than or equal to the largest out-degree in the network. \end{corollary}

Finally, we observe that some of the results of Section~\ref{basicresults}  can be obtained as corollaries of Theorem~\ref{thm:nec_suf_condition}. However, in section~\ref{basicresults} the derivation of these results are  more intuitive and add insight into the relationship between topology and identifiability.

\section{A measurement-based result}\label{measurementbased}
For the application of the results of the previous section, one needs to consider each  node in the network, and examine the paths from its out-neighbors to measured nodes. This procedure may lead one to decide which node needs to be measured. In particular, our results are useful to decide which nodes need to be measured if one wants to identify a particular transfer function. 

In this section, we present a dual approach. We consider a measured node and examine which transfer functions are identifiable from that measured node.
\begin{theorem}\label{measurednode}
Let $j$ be a measured node, and consider a node $i$ that has a directed path to node $j$. Then all transfer functions along that path can be identified if there is no other walk that connects $i$ to $j$.
\end{theorem}
\ni {\bf Proof:}
Let  $N^*_i$ of Theorem~\ref{prop:suf_condition_partial} contain only the out-neighbor of node $i$ that is on the path to $j$ mentioned in the theorem, and let ${\cal C^*}$ contain only $j$. There is no path from any node  $k\in N^+_i\setminus N^*_i$ to $j$, since this path concatenated with the edge $(i,k)$ would constitute another walk from $i$ to $j$. The result then follows from Theorem~\ref{prop:suf_condition_partial}. \cqfd

\begin{figure}[h]
 \centering
\includegraphics[width=0.7\linewidth]{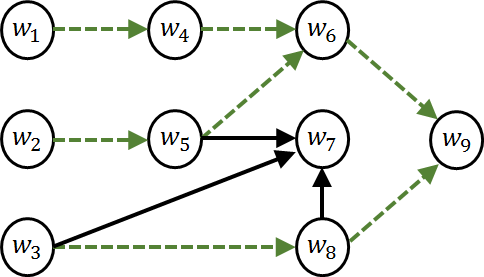}     
\caption{Example illustrating Theorem~\ref{measurednode}: all transfer functions on the dashed green edges can be identified from the measurement of node 9.}
\label{fig:migillustration}
\end{figure}

Theorem~\ref{measurednode} is illustrated by the following example in Figure~\ref{fig:migillustration}; remember again that known signals $r_i$ are added to each node, which are not represented on the figure. It follows from this theorem that the  7 transfer functions on the dashed green-colored  paths can all be identified from the measurement of node 9. If in addition node 7 is also measured, then the 10 transfer functions of the network can all be identified from the two measured nodes 7 and 9.\\

\section{Conclusions}\label{conclusion}

All results so far on the identifiability of networks of dynamical systems have been built on the assumption that all nodes are measured. In this paper, we have addressed the situation where not all nodes are measured,  under the assumptions that the topology is known, that the nodes are noise-free and that they are all excited by known external signals.

We have first shown that the node measurements needed for network identifiability depend entirely on the topology of the network. In doing so, we have observed that the measurement of some nodes (the sinks) are indispensable while other nodes (the sources) yield no information. 
We have then provided a series of results on identifiability. Some of these are based on looking at a particular node and its out-neighbors, and their paths to measured nodes; others have addressed the question of which transfer functions can be identified from the measurement of a particular node.

An important outcome of our work is that networks can often be identified by measuring only a small subset of  nodes. Our main result, based on  the first approach, is a necessary and sufficient condition for identifiability of the transfer functions of all edges leaving a particular node. This result paves the way for an algorithm that lets one decide which transfer functions can be identified from a given set of measured nodes. Our results also yield an easy solution to the problem of identifying just one transfer function embedded in the network. These and other extensions of the present work can be found in a  comprehensive journal paper \cite{Hendrickx&Gevers&Bazanella:17}.

Future research questions will include the search for the sparsest set of measured nodes that allow identification of the whole network as well as the handling of noise.

\section{Appendix}
This section contains the proofs of the results whose proof is not in the main text, as well as some  technical lemmas required for these proofs.\\
\vspace{-2mm}

\ni {\bf Proof of Theorem~\ref{basic1}}\\
1) The first part follows from the definition of a source and from the calculation of $T$ from such $G$ using (\ref{Tdef}). The second part follows  from  (\ref{TzeroT2}). The only way to identify the transfer function on an outgoing path from a source $i$ is if  an external input  signal $r_i$ is applied at the source. \\
2) The first part follows from the definition of a sink and from the calculation of $T$  using (\ref{Tdef}). Let node  $i$ be a sink and let node $k$ be connected to $i$ by a nonzero transfer function $G_{ik}$. Since node $i$ is a terminal node of the path from $k$ to $i$, no  node other than $i$ can give any information about  $G_{ik}$. On the other hand, applying an excitation signal $r_i$ to sink $i$ yields no information, since no path leaves node $i$.\\
\vspace{-2mm}

\ni {\bf Proof of Theorem~\ref{treetheo}}\\
By Theorem~\ref{basic1}  it is necessary to measure all the sinks for any graph, so it remains to prove
sufficiency. Consider an edge $(i,k)$. Since the network contains no cycle, following an arbitrary path form $k$ will eventually lead to some sink $j$. By the assumption of this theorem, there is no other path from $i$ to $j$ than that going through $k$, and neither $k$ nor $i$ are included in any cycle. One can then verify that the closed-loop transfer function $T_{ji}$ from input $i$ to the output $j$ satisfies $T_{ji}=T_{jk}G_{ki}$. Remembering that all transfer functions  $T_{ji}$ and $T_{jk}$ are identifiable (as are all closed-loop transfer function from any input to any sink), we can recover $G_{ki}=T_{ji}/T_{jk}$.
\cqfd 
\vspace{-4mm}

\ni {\bf Proof of Theorem~\ref{theoloop}}\\
Let $\eta$ be the cardinality of $\cal I$ and consider,  without loss of generality,
that the nodes in the cycle are labeled $i=1,\ldots, \eta$ sequentially, that is, there is a link from each node
$i$ to node $i+1$, so that the $\eta$ transfer functions to be identified in the cycle are
$G_{i+1,i}, \; i=1,\ldots , \eta-1$ and $G_{1,\eta}$. Since an external excitation signal is assumed to enter each node, input-output  identification provides all closed-loop transfer functions 
$T_{i,j}, i,j\in {\cal I}$, none of which are zero. Indeed,
\begin{eqnarray*}
&&\hspace{-5mm} T_{i,j} = \frac{1}{\Delta} G_{i,i-1}G_{i-1,i-2} \ldots G_{1,\eta} G_{\eta,\eta-1} \ldots G_{j+1,j}, \; i < j  \\
&&\hspace{-5mm} T_{i,i} = \frac{1}{\Delta} \\
&&\hspace{-5mm} T_{i,j} = \frac{1}{\Delta} G_{i,i-1}G_{i-1,i-2} \ldots  G_{j+1,j}, \; i > j
\end{eqnarray*}
where $\Delta = 1- G_{1,\eta} \Pi_{i=1,\ldots , \eta-1} G_{i+1,i} $. Now, suppose we measure only the
``last" node $i=\eta$. Then  we have identified all the transfer functions
\begin{eqnarray*}
&&\hspace{-5mm} T_{\eta,j} = \frac{1}{\Delta} G_{\eta,\eta-1}G_{\eta-1,\eta-2} \ldots G_{j+1,j}, \; j=1,\ldots,\eta-1 \\
&&\hspace{-5mm} T_{\eta,\eta} = \frac{1}{\Delta} 
\end{eqnarray*}
Now, notice that 
$$
G_{j+1,j} = \frac{T_{\eta,j}}{T_{\eta,j+1}}, \; j=1,\ldots,\eta-1
$$
which gives each one of the $G_{k,l}$ in the loop from $G_{2,1}$ to $G_{\eta,\eta-1}$. The same reasoning holds if we measure  any other node, since it is just a question of relabeling the nodes. \cqfd
For the next lemmas, we use the adjacency matrix $A$ of a directed graph, which we define by $A_{ij}=1$ if there is an edge from $j$ to $i$ and 0 else. 
The following Lemma is standard in graph theory and can easily be proved by recurrence on $k$. 

\begin{lemma}\label{lem:A^k->path}
Let $A$ be the unweighted adjacency matrix of a directed graph. For any integer $k$, $[A^k]_{ij}$ denotes the number of walks of length $k$ from $j$ to $i$. In particular $[A^k]_{ij}=0$ for all $k$ if there is no path from $j$ to $i$.
Similarly, $[G^k]_{ij} = 0$ for all $k$ if there is no path from $j$ to $i$.
\end{lemma}
\vspace{2mm}

\begin{lemma}\label{lem:disjoint_path-->nonsingular}
Let $A$ be the unweighted adjacency matrix of a directed graph consisting of $d$ vertex-disjoint directed paths, respectively from nodes $s_1$ to $t_1$, $s_2$ to $t_2$, \dots $s_d$ to $t_d$.
Then $[(1-A)^{-1}]_{t_is_i} = 1$, and $[(1-A)^{-1}]_{t_is_j} = 0$ if $i\neq j$.
As a consequence, the restriction of $(1-A)^{-1}$ to the $d$ lines  $t_1,t_2,\dots$ and the $d$ rows $s_1,s_2,\dots$ is a permutation matrix, which is nonsingular (determinant 1 or -1).
\end{lemma}
{\bf Proof:}
Let $q$ be the length of the longest of the $d$ disjoint paths in the graph. Since these paths are directed and disjoint, there is no walk of length larger than $q$, hence it follows from Lemma \ref{lem:A^k->path} that $A^{q+1}= 0$. Therefore, 
\begin{equation}\label{eq:inv(I-A)}
(I-A)^{-1} = I + A +A^2 + \dots + A^q
\end{equation}
It then follows from \eqref{eq:inv(I-A)} and Lemma \ref{lem:A^k->path} that $[(I-A)^{-1}]_{ij}$ is the number of walks of length $L$ or less between $j$ and $i$. In particular, (i) there is no walk between $s_i$ and $t_j$ ($j\neq i$), so $[(1-A)^{-1}]_{t_is_j} = 0$, and (ii) there is exactly one walk of length $q$ or less from $s_i$ to $t_i$, so $[(1-A)^{-1}]_{t_is_i} = 1$.
The restriction of $(1-A)$ to the rows $t_1,t_2,\dots$ and columns $s_1,s_2,\dots$ contains only values 0 and 1, with exactly one 1 in each row and each column, and is thus a permutation matrix, which is  nonsingular (Stated otherwise, for each row there is exactly one column with a 1, and the rest is 0).
\cqfd
\vspace{-2mm}

\ni {\bf Proof of Theorem~\ref{prop:suf_condition}}\\
Since our aim is to prove a generic condition, we just need to show that the system in condition \eqref{eq:condition_columnn_detail} has a rank $d^+_i$ for at least one choice of transfer functions.
For the considered node $i$, define $\bar G_{kl} = 1$ if the edge $(l,k)$ belongs to one of the paths mentioned in the claim, and $0$ else. 
The matrix $\bar G$ is then the unweighted adjacency matrix of the subgraph containing only the vertex disjoint paths mentioned in the statement of the proposition, and is also a particular possible instance of $G^0$. Now define $T \bydef (I-\bar G)^{-1}$. 
It then  follows  from Lemma \ref{lem:disjoint_path-->nonsingular} that the restriction $\tilde T$ of $T$ to the rows corresponding to the end points of the paths (which are measured) and the columns corresponding to their origins (i.e. the out-neighbors of $i$) has full rank, i.e. its rank is $d^+_i$. (Note that since the paths from the $d^+_i$ out-neighbors of $i$ are disjoint, they define $d^+_i$ measured nodes, and hence this restriction is a square matrix.)
But this matrix $\tilde T$ is a submatrix of the matrix $T$ in the system
\begin{equation}\label{eq:P_detail}
\sum_{k\in N^+_i} T_{jk}\Delta_{ki} = 0, \forall j\in {\cal C} 
\end{equation}
Thus, this matrix $T$  has full row rank, and hence  condition  (\ref{eq:P_detail}) implies $\Delta_{ki} = 0, \forall k\in N^+_i.$
Since this matrix is obtained by simple (analytic) operations from the $G_{ij}$, having full row rank is a generic property. Hence the fact that $T$ has full row rank for $\bar G$  (i.e. one particular choice of $G^0$ consistent with the network structure), implies it has full row rank for almost all choices of $G^0$ consistent with the network topology. \cqfd
\vspace{-4mm}

\ni {\bf Proof of Theorem~\ref{prop:suf_condition_partial}}\\
It follows from the second part of Lemma \ref{lem:A^k->path} and from $T^0 = (I-G^0)^{-1} = I + G^0 + (G^0)^2 + \dots$ that condition (ii) implies  $T^0_{jk}=0$ if $j\in C^*$ and $k\in N^+_i\setminus N^*_i$. Hence the restriction of the equation system in \eqref{eq:condition_columnn_detail} to the rows of ${\cal C^*}$ does not have any nonzero coefficient for nodes in $N^+_i\setminus N^*_i$, and can be written
\begin{equation}
\sum_{k\in N^*_i} T^0_{jk}\Delta_{ki} = 0, \forall j\in {\cal C^*} 
\end{equation}
The same argument as in the proof of Theorem~\ref{prop:suf_condition} shows that the matrix of this system is generically nonsingular. As a result  $C(I-G)^{-1} =C(I-G^0)^{-1}$ implies   that $\Delta_{ki}=0$ for every $k\in N^*_i$, i.e. we recover all transfer functions on the edges that link $i$ to the nodes of $N^*_i$. \cqfd
\vspace{-4mm}

\ni {\bf Proof of Theorem~\ref{prop:nec_condition}}\\
For a given node $i$, let us partition the indices in three sets $S,B,P$ in the following way: $B$ is the set described in the hypothesis of this theorem. $S$ is the set of nodes that can be reached from a node in $N^+_i$  without going through any node  of $B$ (it does thus not necessarily contain all the nodes of $N^+_i$, as some of them might be in $B$), and $P$ is the set of remaining nodes. By construction of this partition, (i) there is no edge joining $S$ directly to $P$,  (ii) all measured nodes are either in $P$ or $B$, and (iii) all neighbors of $i$ are either in $S$ or $B$. It follows from (ii) and (iii) that the matrix of the system in \eqref{eq:condition_columnn_detail} is a submatrix of $T_{P\cup B, S\cup B}$. And it follows from (i) and Lemma~\ref{prop:partition_matrix} applied to $G$ with the partition $S,B,P$ that the rank of $T_{P\cup B, S\cup B}$ is generically at most $b< |d^+_i|$. Hence condition \eqref{eq:condition_columnn_detail} is generically not satisfied. \cqfd 
\vspace{-4mm}

\begin{lemma}\label{prop:partition_matrix}
Let $G(q)$ be a $L\times L$ transfer function matrix and $\{1,\dots L\} = S \cup B \cup P$ be a partition of the indices such that every path from $S$ to $P$ goes through one of the $b:=|B|$ nodes of $B$,  
i.e. $G_{ij} = 0$ if $i\in P, j\in S$, or in shorthand notation  $G_{PS} = 0$.
Suppose $(I-G)$ is invertible and that $I-G_{PP}$ is invertible (which is generically the case), and let $T = (I-G)^{-1}$ as before. Then the matrix $T_{P \cup B, S\cup B}$ has rank at most $b$.
\end{lemma}
\ni {\bf Proof:}
After re-ordering of the indices, the matrices are 
\beqna
G&=&\left(\begin{array}{ccc}
G_{PP}& G_{PB} & 0\\
G_{BP}& G_{BB} & G_{BS}\\
G_{SP}& G_{SB} & G_{SS}\\
\end{array}
\right)
\text{    and   }\\
T &\bydef& (I-G)^{-1}=\left(\begin{array}{ccc}
T_{PP}& T_{PB} & T_{PS}\\
T_{BP}& T_{BB} & T_{BS}\\
T_{SP}& T_{SB} & T_{SS}\\
\end{array}
\right)
\eeqna
We  focus on the rows $P$ and lines $S$, keeping in mind that $T = I+ GT$. We have
$$
T_{PS} = I_{PS} + [GT]_{PS} = 0 +  G_{PP}T_{PS} + G_{PB} T_{BS} + 0 T_{SS}
$$
Remembering that $(I-G_{PP})$ is assumed invertible, it follows that
\begin{equation}\label{eq:P_TS}
T_{PS} = (I-G_{PP})^{-1} G_{PB} T_{BS}.
\end{equation}
Similarly, there holds
$$
T_{PB} = I_{PB} + [GT]_{PB} = 0 +  G_{PP}T_{PB} + G_{PB} T_{BB} + O T_{SB},
$$
from which follows
\begin{equation}\label{eq:P_TB}
T_{PB} = (I-G_{PP})^{-1} G_{PB} T_{BB}.
\end{equation}
We then obtain  from \eqref{eq:P_TS} and \eqref{eq:P_TB}
\begin{equation}\label{eq:P_TB_BS}
\left( 
\!\!\begin{array}{cc}
T_{PB} & T_{PS}\\
T_{BB} & T_{BS}  
\end{array}
\!\! \right)
=
\left( 
\!\! \begin{array}{c}
(I-G_{PP})^{-1} G_{PB}\\ I
\end{array}
\!\! \right)
\left( 
\begin{array}{cc}
T_{BB} & T_{BS}
\end{array}
\!\! \right),
\end{equation}
which proves the claim of the Lemma since $(T_{BB}, T_{BS})$ has a rank at most $b$, its number of rows.

\end{document}